\documentclass{article}
\usepackage{lscape}

\def\be{\begin{equation}}
\def\ee{\end{equation}}
\def\pab#1#2{\frac{\partial #1}{\partial #2}}
\newcommand{\brak}[1]{\left\{ \begin{array}{lllll} #1 \end{array} \right. }
\newcommand{\ff}[1]{{\mbox{\boldmath $#1$}}}
\def\pab#1#2{\frac{\partial #1}{\partial #2}}
\def\lam{\lambda}

\def\x{\ff{x}}

\begin{document}

\title{True Global Optimality of the Pressure Vessel Design Problem:
A Benchmark for Bio-Inspired Optimisation Algorithms }

\author{Xin-She Yang, Christian Huyck, Mehmet Karamanoglu, Nawaz Khan \\
School of Science and Technology, Middlesex University, \\
The Burroughs, London NW4 4BT, UK. }

\date{}

\maketitle
%% Begin of Main Text %%

\begin{abstract}
The pressure vessel design problem is a well-known design benchmark
for validating bio-inspired optimization algorithms.
However, its global optimality is not clear and
there has been no mathematical proof put forward. In this paper, a detailed
mathematical analysis of this problem is provided that proves that
6059.714335048436 is the global minimum. The Lagrange multiplier
method is also used as an alternative proof and this method is
extended to find the global optimum of a cantilever beam design
problem.
\end{abstract}

{\bf Citation detail:} Xin-She Yang, Christian Huyck, Mehmet Karamanoglu, Nawaz Khan,
rue Global Optimality of the Pressure Vessel Design Problem:
A Benchmark for Bio-Inspired Optimisation Algorithms, 
{\it Int. J. Bio-Inspired Computation}, vol. 5, no. 6, pp. 329-335 (2013).

\section{Introduction}

Engineering optimization is often non-linear with complex constraints,
which can be very challenging to solve. Sometimes, seemingly simple
design problems may in fact be very difficult indeed. Even in very simple
cases, analytical solutions are usually not available, and researchers
have struggled to find the best possible solutions. For example, the
well-known design benchmark of a pressure vessel has only four design
variables [\cite{Cag,GandomiYang,Yang}]; however, the global optimum solution
for pressure vessel design benchmark is still unknown to the research community,
despite a large number of attempts and studies, i.e. \cite{Annaratone,Deb}.
Thus, a mathematical analysis will help to gain some
insight into the problem and thus guide researchers to validate if
their solutions are globally optimal.  This paper attempts to provide
a detailed mathematical analysis of the pressure vessel problem and
find its global optimum. To the best of the author's knowledge, this is a novel
result in the literature.

Therefore, the rest of the paper is organised as follows: Section 2,
introduces the basic formulation of the pressure vessel design
benchmark and then highlights the relevant numerical results from the
literature.  Section 3 provides an analysis of the global
optimum for the problem, whereas Section 4 uses Lagrange
multipliers as an alternative method to prove that the
analysis in Section 3 indeed gives the global optimum. Section
5 extends the same methodology to analyse the optimal solution of
another design benchmark: cantilever beam. Finally, we draw our
conclusions in Section 6.

\section{Pressure Vessel Design Benchmark}

Bio-inspired optimization algorithms have become popular,
and many new algorithms have emerged in recent years
[\cite{YangBA,CheCui,Cui,Gandomi2012ES,YangMOCS,YangMOFA}].
In order to validate new algorithms,  a diverse set of
test functions and benchmarks are often used \cite{Gandomi2011Yang,Jamil}.
Among the structural design benchmarks, the pressure vessel
design problem is one of the most widely used.

In fact, the pressure vessel design problem is a well-known benchmark for
validating optimization algorithms [\cite{Cag,Yang}].  It has four
design variables: thickness ($d_1$), thickness of the heads ($d_2$),
the inner radius ($r$) and the length ($L$) of the cylindrical
section. The main objective is to minimize the overall cost, under
the nonlinear constraints of stresses and yield criteria. The thickness
can only take integer multiples of $0.0625$ inches.

This optimization
problem can be written as
\[ \textrm{minimize} \;\; f(\x) = 0.6224 d_1 r L + 1.7781 d_2 r^2 \]
\be \qquad  + 3.1661 d_1^2 L + 19.84 d_1^2 r,
\label{eq-obj} \ee

subject to
\be \brak{  g_1(\x) = -d_1 + 0.0193 r \le 0 \\ \\
 g_2(\x) = -d_2 + 0.00954 r \le 0 \\  \\
 g_3(\x) = - \pi r^2 L -\frac{4 \pi}{3} r^3 + 1296000 \le 0 \\ \\
 g_4(\x) =L -240 \le 0. } \ee
\noindent The simple bounds are
\be 0.0625 \le d_1, d_2 \le 99 \times 0.0625, \;\; 10.0 \le r, L \le 200. \ee

\begin{landscape}
\begin{table}
\caption{Summary of main results.
[Results marked with $*$ are not valid, see text].} \label{Table-100}
\begin{tabular}{|l|l||l|l|}
\hline
Authors & Results & Authors & Results \\
\hline
\cite{Coello} & 6288.745 &  \cite{LeeGeem} & 7198.433  \\
\cite{Li} & 7127.3 &  \cite{Cai} & 7006.931 \\
\cite{Cag} & {\bf 6059.714} &  \cite{LiCh} & 7127.3 \\
 \cite{CaoW} & 7108.616 &   \cite{Hu} & 6059.131$^*$ \\
 \cite{He} & {\bf 6059.714} &  \cite{CoelloM} & 6059.946 \\
 \cite{Huang} & 6059.734 & \cite{HeW} & 6061.078 \\
 \cite{Liti} & 7197.7 &  \cite{Coello2} & 6263.793 \\
 \cite{Sandgren} & 7980.894 & \cite{Kannan} & 7198.042 \\
 \cite{Akhtar} & 6171 & Yun \cite{Yun} &  7198.424 \\
\cite{Tsai} & 7079.037 &  \cite{Cao} & 7108.616 \\
\cite{Deb} & 6410.381 &  \cite{Coello99} & 6228.744 \\
\cite{Montes} & 6059.702$^*$ &  \cite{Parso} & 6544.27 \\
\cite{Shih} & 7462.1 &   \cite{Kaveh} & 6059.73 \\
\cite{Sandgren2} & 8129.104 & \cite{Coelho} & {\bf 6059.714} \\
 \cite{Wu} & 7207.494 &   \cite{Ray} & 6171 \\
 \cite{Zhang} & 7197.7 &  \cite{CoelloCort} & 6061.123 \\
\cite{Joines} & 6273.28 &  \cite{Michal} & 6572.62 \\
 \cite{Hadj} &  6303.5 &  \cite{Fu} & 8048.6 \\
 \cite{YangBA} & {\bf 6059.714} &  \cite{GandomiYang} & {\bf 6059.714} \\
\hline
\end{tabular}
\end{table}
\end{landscape}

This is a mixed-integer problem, which is usually challenging to
solve. However, there are extensive studies in the literature, and details
can be found in several good survey papers [\cite{Thanedar,Gandomi2011Yang}].
The main results are summarized in Table \ref{Table-100}. It is worth pointing out
that some results are not valid and marked with * in the footnote. These
seemingly lower results actually violated some constraints and/or used
different limits.

As it can be seen
from this table, the results vary significantly from the highest value
of 8129.104 by Sandgren \cite{Sandgren} to the lowest value of
6059.714 by a few researchers [\cite{Cag,Coelho,He,GandomiYang,Gandomi2011mixed,YangBA}].
However, nobody is sure that 6059.714 is the globally optimal solution
for this problem.

The best solution by \cite{GandomiYang} and \cite{YangBA} is
\be f_*=6059.714, \ee
with
\be \x_*=(0.8125, \; 0.4375, \; 42.0984, \; 176.6366). \ee

The rest of the paper analyses this problem mathematically and proves
that this solution is indeed near the global optimum and concludes
that the true globally minimal solution is
$f_{\min}=6059.714335048436$ at
\[ \x_*=(0.8125, \; 0.4375, \]
\be \qquad 42.0984455958549, \; 176.6365958424394). \ee

\section{Analysis of Global Optimality}

As all the design variables must have positive values and $f$ is
monotonic in all variables, the minimization of $f$ requires the
minimization of all the variables if there is no constraint.  As there
are 4 constraints, some of the constraints may become tight or
equalities.  As the range of $L$ is $10\le L \le 200$, the constraint
$g_4$ automatically satisfies $L \le 240$ and thus becomes redundant,
which means that the upper bound
for $L$ is
\be L \le 200. \ee

This is a mixed integer programming problem, which often requires
special techniques to deal with the integer constraints. However, as
the number of combinations of $d_1$ and $d_2$ is not huge (just $100^2
=10,000$), it is possible to go through all the cases for $d_1$ and
$d_2$, and then focus on solving the optimization problems in terms of
$r$ and $L$.

The first two constraints are about stresses. In order to satisfy
these conditions, the hoop stresses $d_1/r$ and $d_2/r$ should be as
small as possible. This means that $r$ should be reasonably large.
For any given $d_1$ and $d_2$, the first two constraints
become
\be r \le \frac{d_1}{0.0193}, \quad r \le \frac{d_2}{0.00954}. \ee

So the upper bound or limit for $r$ becomes
\be U_r=\min \{\frac{d_1}{0.0193}, \frac{d_2}{0.00954} \}. \ee

The above argument that $r$ should be moderately high, may imply that
one of the first two constraints can become tight, or an equality.

The third constraint $g_3$ can be rewritten
\be \pi r^2 L + \frac{4 \pi}{3} r^3 \ge K, \quad K=1296000.
\label {eq-constraint3}\ee
In fact, this is essentially the requirement that the volume of the
pressure vessel must be greater than a fixed volume.
This provides the lower boundary in the search domain of $(r,L)$.

Since $f(r,L)$ is monotonic in $r$ and $L$, the global solution must
be on the lower boundary for any given $d_1$ and $d_2$. In other
words, the inequality $g_3$ becomes an equality
\be \pi r^2 L +\frac{4 \pi}{3} r^3 = K.    \label{equ-g3} \ee

%undone how did this get done?  It works with equation 9, but
%what was the method.

Using equation (\ref {eq-constraint3}), with $L \le 200$, $r$ can be derived
using Newton's method (or an online polynomial root calculator).  The
only positive root implies that
\be r \ge 40.31961872409872=r_{1}. \ee
Similarly, $L \ge 10$ means that
\be r \le 65.22523261350128=r_{2}. \ee
So the true value of $r$ must lie in the interval of $[r_1, r_2]$.

From the first inequality with $r=r_1$, we have
\be d_1 \ge 0.7782. \ee
The second inequality gives
\be d_2 \ge 0.3846. \ee
As both $d_1$ and $d_2$ must be integer multiples $I$ and $J$, respectively,
of $d=0.0625$, the above two inequalities mean
\be I=\lceil \frac{0.7782}{0.0625} \rceil =13, \quad J=\lceil \frac{0.3846}{0.0625} \rceil=7. \ee
In other words, we have
\be d_1 \ge 13 d=0.8125, \quad d_2 \ge 7 d=0.4375. \ee

From the objective function (Eq. \ref {eq-obj}), both $d_1$ and $d_2$
should be as small as possible, so as to get the minimum possible
$f$. This means that the global minimum will occur at $d_1=0.8125$ and
$d_2=0.4375$.

Now the objective function with these $d_1$ and $d_2$ values can be written as
\[ f(r,L)=0.5057 r L + 0.77791875 r^2 \]
\be \qquad + 2.090120703125 L + 13.0975 r. \label{new-obj} \ee

As $d_1=0.8125$ and $d_2=0.4375$, the first inequalities ($g_1$ and
$g_2$) will give an upper bound of $r$
\[ R_*=\min \{ \frac{d_1}{0.0193}, \frac{d_2}{0.00954} \} \]
\[ \qquad =\min \{42.0984455958549, 45.859538784067 \} \]
\be \qquad =42.0984455958549. \ee

Again from the objective function, which is monotonic in terms of $r$ and $L$,
the optimal solution should  occur at the two extreme ends of the boundary governed by Eq.~(\ref{equ-g3}).

The one end at $r=R_*$ gives
\be L_*=\frac{K}{\pi R_*^2}-\frac{4 R_*}{3}=176.6365958424394. \ee
This is the point for the global optimum with
\be f_{\min}=6059.714335048436. \label{fmin-equ}  \ee
The other extreme point is at $r'=40.31961872409872$ and $L=200$, which leads to an objective value of
\be f'=6288.67704565344, \ee
and clearly is not the global optimum.

\section{Method of Lagrange Multipliers}

The optimal solution (\ref{fmin-equ}) can alternatively be proved by
solving the following constrained problem with one equality because
all of the upper bounds or inequalities are automatical satisfied.  By minimizing
$d_1$, and $d_2$, the objective function becomes:
\[ \textrm{minimize }\;\;  f(r,L)=0.5057 r L + 0.77791875 r^2 \]
\be \qquad + 2.090120703125 L + 13.0975 r. \label{new-obj2} \ee
subject to
\be g_e = \pi r^2 L + \frac{4 \pi}{3} r^3 -K =0, \ee
with the simple bounds
\[ 40.31961872409872 \le r \le 42.098445595854919, \]
\be \qquad 10 \le L \le 176.6365958424394. \ee
To avoid writing long numbers, let us define
\[ a=0.5057, \quad b= 0.77791875, \]
\be \qquad c=2.090120703125, \quad d=13.0975. \ee
We have
\be \textrm{minimize }\;\;  f(r,L)=a r L + b r^2 + c L + d r. \label{new-obj2} \ee

This problem can be solved by the Lagrange multiplier method, and we have
\be \textrm{minimize } \phi=f+\lam g_e, \ee
where $\lam$ is the Lagrange multiplier.

The optimum should occur when
\[ \pab{\phi}{r} = \pab{f}{r}+\lam \pab{g_e}{r} \]
\be \quad = a L + 2b  r + d + \lam (2 \pi r L + 4 \pi r^2) =0, \ee

\be \pab{\phi}{L} = \pab{f}{L} +\lam \pab{g_e}{L} = a r + c  + \lam (\pi r^2)=0, \ee
\be \pab{\phi}{\lam}=\pi r^2 L +\frac{4 \pi}{3} r^3 - K=0. \ee

Now  we have three equations for three unknowns
\be \brak{a L + 2b r + d + \lam (2 \pi r L + 4 \pi r^2)=0, \\ \\
a r + c + \lam (\pi r^2)=0, \\ \\
\pi r^2 L +\frac{4 \pi r^3}{3} -K=0.
} \ee

The equation in the middle gives
\be \lam =- \frac{a r + c}{\pi r^2}, \ee

Substituting this, together with $L=K/(2r^2) -4 r/3$, into the first equation, we have
\[ a (K r -\frac{4 \pi r^3}{3}) + 2 b \pi r^4 + \pi d r^3 \]
\be \qquad -(a r+c) (2 K +\frac{4 \pi r^3}{3})=0, \ee
which is a quartic equation with four roots in general. The only feasible solution
within $[r_1, r_2]$ is $42.098445595854919$, which corresponds to
$L=176.6365958424394$. This solution is indeed the global best solution as given in (\ref{fmin-equ}).

\section{Cantilever Beam Design Benchmark}

Another widely used benchmark for validating bio-inspired algorithms is the
design optimization of a cantilever beam, which is to minimize the overall weight of a
cantilever beam with square cross sections [\cite{Fleu,GandomiYang,Gandomi2011mixed}].
It can be formulated as
\be \textrm{minimize }  f(\x)=0.0624 (x_1+x_2+x_3+x_4+x_5), \ee
subject to the inequality
\be g(\x)=\frac{61}{x_1^3}+\frac{37}{x_2^3}+\frac{19}{x_3^3} + \frac{7}{x_4^3}+\frac{1}{x_5^3}-1 \le 0. \ee
The simple bounds/limits for the five design variables are
\be 0.01 \le x_i \le 100, \quad i=1,2,...,5. \ee

Since the objective $f(\x)$ is linear in terms of all design
variables, and $g(\x)$ encloses a hypervolume, it can be thus expected
that the global optimum occurs when the inequality becomes
tight. That is, the inequality becomes an equality
\be g(\x)=\frac{61}{x_1^3}+\frac{37}{x_2^3}+\frac{19}{x_3^3} + \frac{7}{x_4^3}+\frac{1}{x_5^3}-1 =0. \ee

For ease of analysis, we rewrite the above equation as
\be g(\x)=\sum_{i=1}^5 \frac{a_i}{x_i^3}-1=0, \ee
where
\be {\bf a}=(a_1, a_2, a_3, a_4, a_5)=(61, 37, 19, 7, 1). \ee
Hence, the cantilever beam problem becomes
\be \textrm{minimize } \;\; f(\x)=k \sum_{i=1}^5 x_i, \quad k=0.0624, \ee
subject to
\be g(x)=\sum_{i=1}^5 \frac{a_i}{x_i^3}-1=0. \ee
Using the method of Lagrange multipliers, we have
\[ \textrm{Minimize } \phi=f+\lam g \]
\be \qquad =k \sum_{i=1}^5 x_i +\lam \Big( \sum_{i=1}^n \frac{a_i}{x_i^5}-1 \Big). \ee
Then, the optimality conditions give
\be \pab{\phi}{x_i} = k + \lam (-3) \frac{a_i}{x_i^4}=0, \quad i=1,2,...,5, \label{opt-100} \ee
\be \pab{\phi}{\lam}=\sum_{i=1}^5 \frac{a_i}{x_i^3}-1=0. \label{opt-200} \ee
From Eq. (\ref{opt-100}), we have
\be x_i^4 =\frac{3 \lam a_i}{k}, \quad \textrm{ or } \quad \frac{a_i}{x_i^3} =\frac{k x_i}{3 \lam}. \label{opt-300} \ee
Substituting it into Eq. (\ref{opt-200}), we have
\be \sum_{i=1}^5 (\frac{k x_i}{3 \lam}) -1 = \frac{k}{3 \lam} (\sum_{i=1}^5 x_i)-1=0. \ee
After rearranging and using results from (\ref{opt-300}), now we have
\be \frac{3 \lam}{k}=\sum_{i=1}^5 \Big(\frac{ 3 \lam a_i}{k}\Big)^{1/4},  \ee
which is a nonlinear equation for $\lam$.
However, it is straightforward to find that
\be \lam \approx 0.4466521202, \ee
which leads to the optimal solution
\[ \x_*=(6.0160159, 5.3091739, \]
\be \qquad 4.4943296, 3.5014750, 2.15266533), \ee
with the minimum
\be f_{\min}(\x_*)=1.339956367. \ee
This is the global optimum. However, the authors have not seen any studies that have found this solution in the literature. Slightly higher values have been found by cuckoo search and other methods
[\cite{Chicke,GandomiYang}]. The best solution found so far by \cite{GandomiYang} is
\be \x_{\rm best}=(6.0089, 5.3049, 4.5023, 3.5077, 2.1504), \ee
and
\be f_{\rm best}=1.33999, \ee
which is near this global optimum.

The above mathematical analysis can be very useful to guide future validation of new optimization methods when the above design benchmarks are used.

\section{Conclusions}

Pressure vessel design problem is a well-tested benchmark that has been used for validating
optimization algorithms and their performance. We have provided a detailed mathematical analysis and obtained
its global optimality. We have also used the method of Lagrange multipliers to double-check that
the obtained optimum is indeed the global optimum for the pressure vessel design problem.
By using the same methodology, we also analysed the design optimization of a cantilever beam.

However, it is worth pointing out that the method of Lagrange multipliers is only valid for
optimization problems with equalities or when an inequality becomes tight. For general nonlinear
optimization problems, we have to use the full Karush-Kuhn-Tucker (KKT) conditions to analyze their
optimality [\cite{Yang}], though such KKT can be extremely challenging to analyse in practice.

Even for design problems with only a few design variables, an analytical solution will provide
greater insight into the problem and thus can act as better benchmarks for validating new
optimization algorithms. Further work can focus on the analysis of other nonlinear design
benchmarks.


\begin{thebibliography}{10}



\bibitem[Akhtar et al.~(2002)]{Akhtar}
Akhtar, S., Tai, K., Ray, T. (2002). A socio-behavioural simulation
model for engineering design optimization, {it Engineering Optmization},
{\bf 34}(4), pp. 341--354.

\bibitem[Annaratone~(2007)]{Annaratone}
Annaratone, D. (2007) Pressure Vessel Design, Springler-Verlag Berlin Heidelberg.


\bibitem[Cagnina et al.~(2008)]{Cag}
Cagnina L. C., Esquivel S. C., and Coello C. A., (2008).
Solving engineering optimization problems with the simple
constrained particle swarm optimizer, {\it Informatica},
{\bf 32}, 319-326.

\bibitem[Cai and Thierauf~(1997)]{Cai}
Cai, J., Thierauf, G., (1997). Evolution strategies in engineering
optimization, {\it Eng Optimization}, {\bf 29}(1), pp. 177--199.

\bibitem[Cao and Wu~(1997)]{Cao}
Cao, Y. J., Wu, Q. H., (1997). Mechanical design optimization by
mixed variable evolutionary programming. In: {\it Proceedings of the
1997 International Conference on Evolutionary Computation},
Indianapolis, pp. 443--446.

\bibitem[Cao and Wu~(1999)]{CaoW}
Cao, Y. J., Wu, Q. H., (1999). A mixed variable evolutionary programming
for optimization of mechanical design, {\it Int J Eng Intel
Syst Elect Eng Commun},  {\bf 7}(2), pp. 77--82.

\bibitem[Che and Cui~(2011)]{CheCui}
Che, Z. H., Cui, Z. H., (2011).
Unbalanced supply chain design using the analytic network process and a hybrid heuristic-based algorithm with balance modulating mechanism,
{\it Int. J. Bio-inspired Computation}, {\bf 3}(1), 56--66.

\bibitem[Chickermane and Gea~(1996)]{Chicke}
Chickermane, H., Gea, H. C., (1996). Structural optimization using a
new local approximation method, {\it  Int J Numer Method Eng},
{\bf 39}, pp.829--846.


\bibitem[Santos Ceolho~(2010)]{Coelho}
dos Santos Coelho, L, (2010). Gaussian quantum-behaved particle
swarm optimization approaches for constrained engineering
design problems, {\it  Expert. Syst. Appl.}, {\bf 37}(2), pp. 1676�1683.


\bibitem[Coello~(2000a)]{Coello}
Coello C. A. C., (2000a).  Use of a self-adaptive penalty approach for
engineering optimization problems. {\it Comput. Ind.}, {\bf 41}(2), pp. 113-�127.

\bibitem[Coello~(2000b)]{Coello2}
Coello, C. A. C., (2000b). Constraint-handling using an evolutionary
multiobjective optimization technique, {\it Civil Engrg Environ Syst},
{\bf 17}, pp. 319--346.

\bibitem[Coello~(1999)]{Coello99}
Coello C. A. C., (1999).  Self-adaptive penalties for GA based optimization,
{\it Proc. Congr. Evol. Comput.}, {\bf 1}, pp. 573--580.

\bibitem[Coello and Cort\'es~(2004)]{CoelloCort}
Coello, C. A. C.,  Cort\'es, N. C., (2004). Hybridizing a genetic algorithm
with an artificial immune system for global optimization, {\it Engineering
Optmization}, {\bf 36}(5), pp. 607�634.

\bibitem[Coello and Mezura Montes~(2001)]{CoelloM}
Coello, C. A. C., Mezura Montes, E., (2001). Use of dominance-based
tournament selection to handle constraints in genetic algorithms. In:
{\it Intelligent Engineering Systems through Artificial Neural Networks}
(ANNIE2001), {\bf 11}, ASME Press, St. Louis, pp. 177--182.

\bibitem[Cui et al.~(2013)]{Cui}
Cui, Z. H., Fan S. J., Zeng J. C., Shi, Z. Z., (2013). Artificial plant
optimisation algorithm with three-period photosynthesis,
{\it Int. J. Bio-Inspired Computation}, {\bf 5}(2), 133--139.




\bibitem[Deb and Gene~(1997)]{Deb}
Deb, K., Gene, A. S., (1997). A robust optimal design technique for
mechanical component design. in: Evolutionary algorithms in engineering
applications. Springer-Verlag, Berlin, pp. 497--514.

\bibitem[Fleury and Braibant~(1986)]{Fleu}
Fleury, C., Braibant, V., (1986). Structural optimization: a new dual
method using mixed variables, {\it Int J Numer Meth Eng},  {\bf 23}, pp. 409--428.


\bibitem[Fu et al.~(1991)]{Fu}
Fu, J., Fenton, R. G., Cleghorn, W. L., (1991). A mixed integer-discrete continuous
programming method and its application to engineering
design optimization, {\it Engeering Optmization}, {\bf 17}, pp. 263--280.

\bibitem[Gandomi et al.~(2013)]{GandomiYang}
Gandomi, A. H., Yang, X.-S., and Alavi, A. H., (2013).
Cuckoo search algorithm: a metaheuristic approach to solve
structural optimization problems, {\it Engineering with Computers},
{\bf 29}(1), pp. 17-35.

\bibitem[Gandomi et al.~(2011)]{Gandomi2011mixed}
Gandomi, A. H., Yang, X. S., Alavi, A. H., (2011).
Mixed variable structural optimization using Firefly Algorithm,
{\it Computers \& Structures}, {\bf 89} (23-24), 2325--2336.

\bibitem[Gandomi and Yang~(2011)]{Gandomi2011Yang}
Gandomi, A. H. and Yang, X. S., (2011).
Benchmark problems in structural optimization, in: {\it Computational Optimization,
Methods and Algorithms}, Studies in Computational Intelligence (Eds. Koziel S. and Yang X. S.), Springer, Heidelberg, vol. 356, pp. 259--281.

\bibitem[Gandomi et al.~(2012)]{Gandomi2012ES}
Gandomi, A. H., Yang, X. S., Talatahari, S., Deb, S., (2012).
Coupled eagle strategy and differential evolution for unconstrained
and constrained global optimization, {\it Computers \& Mathematics with
Applications}, {\bf 63}(1), 191--200.

\bibitem[Hadj-Alouane and Bean ~(1997)]{Hadj}
Hadj-Alouane, A. B., Bean, J. C., (1997).
A genetic algorithm for the multiple-choice
integer program. {\it Oper Res},  {\bf 45}, pp. 92--101.


\bibitem[He et al.~(2004)]{He}
He, S., Prempain, E., Wu, Q. H., (2004). An improved particle swarm
optimizer for mechanical design optimization problems, {\it Engineering
Optimization}, {\bf 36}(5), pp. 585--605.

\bibitem[He and Wang~(2006)]{HeW}
He, Q., Wang, L., (2006). An effective co-evolutionary particle
swarm optimization for engineering optimization problems, {\it Eng
Appl Artif Intel},  {\bf 20}, pp.89--99.


\bibitem[Hu et al.~(2003)]{Hu}
Hu, X., Eberhart, R. C., Shi, Y. (2003). Engineering optimization with
particle swarm. In: {\it Proc. 2003 IEEE Swarm Intelligence Symposium},
pp. 53--57.

\bibitem[Huang et al.~(2007)]{Huang}
Huang, F. Z., Wang, L., He, Q. (2007). An effective co-evolutionary
differential evolution for constrained optimization, {\it Appl. Math.
Comput.}, {\bf 186}, pp. 340--356.

\bibitem[Jamil and Yang~(2013)]{Jamil}
Jamil, M., and Yang, X. S., (2013). A litrature survey of benchmar functions for
global optimisation problems, {\it Int. J. of Mathematical Modelling and Numerical Optimisation}, {\bf 4}(2), 150--194.




\bibitem[Joines and Houck~(1994)]{Joines}
Joines, J., Houck, C. (1994). On the use of non-stationary penalty
functions to solve nonlinear constrained optimization problems
with GAs. In: {\it Proceedings of the first IEEE Conference on
Evolutionary Computation}, Orlando, Florida. D. Fogel (ed.). IEEE Press, pp 579--584


\bibitem[Kannan and Kramer~(1994)]{Kannan}
Kannan, B. K., Kramer, S. N., (1994). An augmented Lagrange multiplier
based method for mixed integer discrete continuous optimization
and its applications to mechanical design. {\it J Mech Des. Trans.},
{\bf  116}, pp. 318--320

\bibitem[Kaveh and Talatahari~(2010)]{Kaveh}
Kaveh, A., Talatahari, S., (2010). An improved ant colony optimization
for constrained engineering design problems, {\it Engineering Computations},
{\bf 27}(1), pp.155--182.


\bibitem[Lee and Geem~(2005)]{LeeGeem}
Lee K. S., Geem Z. W. (2005), A new meta-heuristic algorithm for
continuous engineering optimization: harmony search theory and
practice, {\it Comput. Methods Appl. Mech. Eng.}, {\bf 194}, pp. 3902--3933.

\bibitem[Li and Chou~(1994)]{Li}
Li,  H.-L., Chou C.-T., (1994). A global approach for nonlinear mixed
discrete programming in design optimization. Engineering Optimization,
{\bf  22}, pp. 109--122.

\bibitem[Li and Chang~(1998)]{LiCh}
Li, H. L., Chang, C. T., (1998). An approximate approach of global
optimization for polynomial programming problems, {\it Eur J. Oper
Res}, {\bf 107}(3), pp.625--632.


\bibitem[Litinetskiand Abramzon~(1998)]{Liti}
Litinetski, V. V., Abramzon, B. M., (1998). Multistart adaptive
random search method for global constrained optimization in
engineering applications, {\it Engineering Optimization}, {\bf 30}(2), pp. 125--154.

\bibitem[Michalewicz and Attia~(1994)]{Michal}
Michalewicz, Z., Attia, N. (1994). Evolutionary optimization of
constrained problems. Proceedings of the 3rd Annual Conference
on Evolutionary Programming, World Scientific, pp. 98--108.

\bibitem[Montes et al.~(2007)]{Montes}
Montes, E. M., Coello, C. A. C., Vel\'azquez-Reyes, J, Mu\~noz-D\'avila, L.,
(2007). Multiple trial vectors in differential evolution for engineering
design, {\it Engineering Optimization}, {\bf 39}(5), pp. 567--589.


\bibitem[Parsopoulos and Vrahatis~(2005)]{Parso}
Parsopoulos,  K. E., Vrahatis,  M. N.,. (2005). Unified particle swarm
optimization for solving constrained engineering optimization
problems. In: {\it Lecture Notes in Computer Science (LNCS)},
{\bf 3612}, pp. 582--591.

\bibitem[Rat and Liew~(2003)]{Ray}
Ray, T., Liew, K., (2003). Society and civilization: An optimization
algorithm based on the simulation of social behavior, {\it IEEE Trans
Evol Comput}, {\bf 7}(4), pp.386--396.

\bibitem[Sandgren~(1998)]{Sandgren}
Sandgren, E., (1988). Nonlinear integer and discrete programming
in mechanical design, {\it Proceedings of the ASME Design Technology
Conference}, Kissimine, FL, pp. 95--105.


\bibitem[Sandgren~(1990)]{Sandgren2}
Sandgren, E. (1990). Nonlinear integer and discrete programming
in mechanical design optimization, {\it  J Mech Design}, {\bf 112}(2), pp. 223--229.

\bibitem[Shih and Lai~(1995)]{Shih}
Shih, C. J., Lai, T. K., (1995). Mixed-discrete fuzzy programming for
nonlinear engineering optimization, {\it Engineering Optimization},
{\bf 23}(3), pp. 187--199.


\bibitem[Thanedar and Vanderplaats~(1995)]{Thanedar}
Thanedar, P. B.,  Vanderplaats, G. N., (1995).
Survey of discrete variable optimization for structural design,
{\it Journal of Structural Engineering ASCE}, {\bf 121} (2), 301--306.


\bibitem[Tsai et al.~(2002)]{Tsai}
Tsai, J.-F., Li, H.-L., Hu, N.-Z., (2002).
Global optimization for signomial discrete programming problems
in engineering design, {\it Engineering Optmization}, {\bf 34}(6), pp. 613--622.

\bibitem[Wu and Chow~(1995)]{Wu}
Wu, S. J., Chow, P. T., (1995). Genetic algorithms for nonlinear mixed
discrete-integer optimization problems via meta-genetic parameter
optimization,  {\it Engineering Optimization}, {\bf 24}, pp. 137--159.

\bibitem[Yang~(2010)]{Yang}
 Yang X. S., (2010). {\it Engineering Optimisation: An Introduction with
 Metaheuristic Applications}, John Wiley and Sons.

\bibitem[Yang and Gandomi~(2012)]{YangBA}
Yang, X. S. and Gandomi, A. H., (2012). Bat algorithm: a novel apporach for
global engineering optimization, {\it Engineering Computations},
{\bf 29}(5), pp. 464--483.

\bibitem[Yang~(2013)]{YangMOFA}
Yang, X. S., (2013).
Multiobjective firefly algorithm for continuous optimization,
{\it Engineering with Computers}, {\bf 29}(2), 175--184.

\bibitem[Yang and Deb~(2013)]{YangMOCS}
Yang, X. S. and Deb, S., (2013).  Multiobjective cuckoo search
for design optimization, {\it Computers \& Operations Research},
{\bf 40}(6), 1616--1624.

\bibitem[Yun~(2005)]{Yun}
Yun,  Y. S., (2005). Study on Adaptive Hybrid Genetic Algorithm and
Its Applications to Engineering Design Problems, Waseda University, MSc Thesis.

\bibitem[Zhang and Wang~(1993)]{Zhang}
Zhang, C., Wang, H. P., (1993). Mixed-discrete nonlinear optimization
with simulated annealing. {\it Engineering Optmization},  {\bf 17}(3), pp. 263�280.


\end{thebibliography}
\end{document}